
\magnification=1100
\baselineskip=14pt

\def\CHA{C}
\def\DUB{D}
\def\Falk{F}
\def\GU{GU}
\def\MON{M}
\def\ROST{R}
\def\RY{RY}

\def\sgn{\mathop{\rm sgn}}

  \def\FF{{\cal F}}

\def\EEE{{\bf E}}       
      
      \def\PPP{{\bf P}} 
  \def\RRR{{\bf R}}

\def\rect#1#2#3{\raise .1ex\vbox{\hrule height.#3pt
   \hbox{\vrule width.#3pt height#2pt \kern#1pt\vrule width.#3pt}
        \hrule height.#3pt}}

\def\qed{$\hskip 5pt\rect364$} 

\def\ref#1{{\rm [}{\bf #1}{\rm ]}}   
\def\nref#1#2{{\rm [}{\bf #1}{\rm ;\ #2]}}

\def\nin{\noindent}
\def\prf{\nin{\it Proof. }}

\def\today{\ifcase\month\or 
  January\or February\or March\or April\or  
  May\or June\or July\or August\or  
  September\or October\or November\or
  December\fi\space\number\day,\ \number\year}

\outer\def\proclaim#1{\medbreak\noindent\bf\ignorespaces
   #1\unskip.\enspace\sl\ignorespaces}
\outer\def\endproclaim{\par\ifdim\lastskip<\medskipamount\removelastskip
   \penalty 55 \fi\medskip\rm}

\centerline{\bf A Peculiar Proof of the Martingale Convergence Theorem}
\medskip
\centerline{P.J. Fitzsimmons}
\medskip
\centerline{Department of Mathematics}

\centerline{University of California, San Diego}

\centerline{9500 Gilman Drive}

\centerline{La Jolla, CA 92093--0112}

\centerline{\tt pfitzsimmons@ucsd.edu}
\bigskip

\centerline{\bf ABSTRACT}
\smallskip
We prove the Martingale Convergence Theorem by using the work of L.~Dubins and I.~Monroe about embedding a given discrete-time martingale in the sample paths of a Brownian motion.
\bigskip\bigskip

\bigskip

We provide a new (if ridiculous) proof of the Martingale Convergence Theorem.

Let $M=(M_n)_{n=0}^\infty$ be a discrete time martingale. With no real loss of generality we assume  $M_0=0$.
We consider first the case of $L^1$-bounded martingales.

\proclaim{(1) Theorem} (J.L. Doob) Suppose that
$$
K:=\sup_n\EEE|M_n|<\infty.
\leqno(2)
$$
Then
$$
M_\infty:=\lim_n M_n
$$
exists in $\RRR$ a.s., and $M_\infty$ is integrable.
\endproclaim

\prf Because $n\mapsto |M_n|$ is a submartingale, we  deduce from (2) that
$$
\EEE|M_T|\le K\leqno(3)
$$
for each a.s.  finite stopping time $T$. 
We know from the work of  Monroe \ref{\MON} and Dubins \ref{\DUB} (see also \nref{\CHA}{Thm.~10.3}) that  the martingale $M$ can be embedded in  Brownian motion as follows. Let $(B_t)$ be a  1-dimensional Brownian motion with $B_0=0$, defined on a filtered probability space $(\Omega,\FF, (\FF_t), \PPP)$ satisfying the usual conditions. Here $(\FF_t)$ denotes the filtration generated by $B$, suitably completed.
Then there is a sequence $0=T_0\le T_1\le T_2\le\cdots$ of standard $(\FF_t)$-stopping times such that the process $M'$ defined by  $M'_n:=B_{T_n}$, $n=0,1,2,\ldots$ 
has the same law as the original martingale $M$. (In this context, we say that a stopping time $T$ is  ``standard" provided the stopped process $t\mapsto B_{t\wedge T}$ is uniformly integrable. In particular, such a stopping time is finite a.s.)  
To show that $\lim_nM_n$ exists a.s., it therefore suffices to show that $M'$ has the same property.
To take advantage of (3) we use Tanaka's formula (see \nref{\RY}{p.~222}):
$$
|B_t|=\int_0^t\sgn(B_s)\,dB_s+L_t,\qquad t\ge 0,
\leqno(4)
$$
where $L$ is the local time of $B$ at $0$. Consequently, if $T$ is an $(\FF_t)$-stopping time, then
$$
\EEE[L_{T\wedge m}]=\EEE|B_{T\wedge m}|,\qquad \forall m>0.
\leqno(5)
$$
In particular,
$$
\EEE[L_{T_n\wedge m}]=\EEE|B_{T_n\wedge m}|,\qquad n=0,1,2,\ldots,
\leqno(6)
$$
and by standardness of $T_n$ and monotone convergence, after a passage to the limit in (6)  as $m\to\infty$ we obtain
$$
\EEE[L_{T_n}]=\EEE|B_{T_n}| =\EEE|M_n|\le K,\qquad n=0, 1,2,\ldots.
\leqno(7)
$$
Let $T_\infty:=\lim_n T_n$. Because $t\mapsto L_t$ is monotone:
$$
\EEE[L_{T_\infty}]\le K<\infty.
\leqno(8)
$$
But $L_\infty=\infty$ a.s. because of the recurrence of one-dimensional Brownian motion. It follows from (8)  that
$$
\PPP[T_\infty<\infty]=1.
$$
Therefore the limit
$$
M'_\infty:=\lim_nM'_n=\lim_nB_{T_n}=B_{T_\infty}
$$
exists in $\RRR$, a.s. The integrability of $M'_\infty$ follows from Fatou's lemma.
\qed
\bigskip

The same sort of proof can be given for the convergence theorem for non-negative supermartingales. The proof is more elementary  because  Tanaka's formula is not needed. Indeed, let $X=(X_n)_{n=0}^\infty$ be a non-negative supermartingale with $X_0=1$ (without loss of generality). Let the Brownian motion $(B_t)$ be as before.
We know from Gushchin and Urusov \ref{\GU} that $X$ can be embedded in the  geometric Brownian motion $Z_t=\exp(B_t-t/2), t\ge 0$. That is, there is an increasing sequence $0=T_0\le T_1\le\cdots\le T_n\le T_{n+1}\le \infty$ of $(\FF_t)$-stopping times  such that $X'_n:=Z_{T_n}$, $n=0,1,2,\ldots$ has the same distribution as $X$. (N.B.: The stopping times $T_n$ are now allowed to take the value $\infty$. 
By the law of large numbers for Brownian motion, $\lim_{t\to\infty}Z_t =0$ a.s., so we take $Z_\infty:=0$.) 
By monotonicity we have  $T_\infty:=\lim_n T_n\in[0,\infty]$, and then by path continuity $\lim_nX'_n=\lim_nZ_{T_n}=Z_{T_\infty}$ a.s.
It follows that $\lim_n X_n$ exists in $[0,\infty[$ a.s. \qed
\bigskip

It may be noted that because the geometric Brownian motion is a  time change of Brownian motion on $[0,\infty[$ absorbed in state $0$, the embedding result for non-negative supermartingales can be restated as follows: As before, let  $B=(B_{t})_{t\ge 0}$, be a  Brownian motion , but now started at $1$; let $\tau:=\inf\{t>0: B_t=0\}$. If  $X=(X_n)_{n=0}^\infty$ is a non-negative supermartingale with $X_0=1$   then there is a sequence $(T_n)$ of stopping times  of $B$ with $0=T_0\le T_1\le T_2\le\cdots\le T_n\le T_{n+1}\le\cdots\le\tau$ such that $(B_{T_n})_{n=0}^\infty$ has the same law as $X$.
\bigskip

The embedding of Gushchin and Urusov  uses randomized stopping times, but this can be avoided when embedding a supermartingale in absorbing Brownian motion by employing Falkner's refinement \nref{\Falk}{Thm.~1} of Rost's theorem \ref{\ROST} for transient Markov processes.
\bigskip

\bigskip
\centerline{\bf References}
\medskip
\frenchspacing

\itemitem{[\CHA]}
Chacon, R.V.: Potental processes, {\it Trans. Amer. Math. Soc.} {\bf 226} (1977) 39--58.
\smallskip

\itemitem{[\DUB]}
Dubins, L.E.: On a theorem of Skorohod, {\it Ann. Math. Statist.} {\bf 39} (1968)  2094--2097. 
\smallskip

\itemitem{[\Falk]}
Falkner, N.: Stopped distributions for Markov processes in duality, {\it Z. Wahrsch. Verw. Gebiete} {\bf 62} (1983) 43--51. 
\smallskip

\itemitem{[\GU]}
Gushchin, A.A., Urusov, M.A.: Processes that can be embedded in a geometric Brownian motion, {\it Theory Probab. Appl.\/} {\bf 60} (2016) 246--262.
\smallskip

\itemitem{[\MON]}
Monroe, I.: On embedding right continuous martingales in Brownian motion, {\it Ann. Probab.\/} {\bf 43} (1972) 1293--1311.
\smallskip

\itemitem{[\RY]}
Revuz, D., Yor, M.: {\it Continuous martingales and Brownian motion,} (3rd ed.).  Springer-Verlag, Berlin, 1999. 
\smallskip

\itemitem{[\ROST]}
Rost, H.: The stopping distributions of a Markov process, {\it Invent. Math.\/} {\bf 14} (1971), 1--16. 
\smallskip

\end